\documentclass[12pt]{article}
\usepackage{amssymb}

\usepackage[top=2.75cm, bottom=2.75cm, left=3.0cm, right=3.0cm]{geometry}

\newcommand{\be}{\begin{equation}}
\newcommand{\ee}{\end{equation}}
\newcommand{\beqs}{\begin{eqnarray}}
\newcommand{\eeqs}{\end{eqnarray}}

\newcommand{\ie}{{\it i.e.},~}

\def\({\left(}
\def\){\right)}

\def\ni{\noindent}

\def\mxth{\mathsurround=0pt }
\def\xversim#1#2{\lower2.pt\vbox{\baselineskip0pt \lineskip-.5pt
x  \ialign{$\mxth#1\hfil##\hfil$\crcr#2\crcr\sim\crcr}}}

\def\be{\begin{equation}}
\def\ee{\end{equation}}
\def\bea{\begin{eqnarray}}
\def\eea{\end{eqnarray}}

\newcommand{\ft}[2]{{\textstyle\frac{#1}{#2}}}

\DeclareSymbolFont{AMSb}{U}{msb}{m}{n}
\DeclareMathSymbol{\fieldC}{\mathalpha}{AMSb}{"43}
\DeclareMathSymbol{\fieldR}{\mathalpha}{AMSb}{"52}
\DeclareMathSymbol{\fieldZ}{\mathalpha}{AMSb}{"5A}

\begin{document} 

\setcounter{page}{0}
\thispagestyle{empty}

\begin{flushright} \small
YITP--SB--06--06 \\ HUTP--06/A0007 \\ ITP--UU--06/07 \\ SPIN--06/05 \\ ITFA-2006-09 \\
\end{flushright}
\smallskip

\begin{center} \LARGE
 Quaternion-K\"ahler spaces, hyperk\"ahler cones, \\ and the c-map
  \\[6mm] \normalsize
 Martin Ro\v{c}ek$^{1,2}$, Cumrun Vafa$^3$ and Stefan Vandoren$^4$ \\[3mm]
 {\small\slshape
 $^1$C.N. Yang Institute for Theoretical Physics\\
 State University of New York at Stony Brook, NY 11790, USA\\[3mm]
 $^2$Institute for Theoretical Physics\\
 University of Amsterdam, Valckenierstraat 65, 1018 XE Amsterdam, The Netherlands\\[3mm]
 $^3$Jefferson Physical Laboratory, Harvard University\\
 Cambridge, MA 02138, USA \\[3mm]
 $^4$Institute for Theoretical Physics \emph{and} Spinoza Institute \\
 Utrecht University, 3508 TD Utrecht, The Netherlands}
\end{center}
\vspace{10mm}

\centerline{\bfseries Abstract} \medskip

\ni Under the action of the c-map, special K\"ahler manifolds 
are mapped into a class of quaternion-K\"ahler spaces. 
We explicitly construct the
corresponding Swann bundle or hyperk\"ahler cone, and
determine the hyperk\"ahler potential in terms of the
prepotential of the special K\"ahler geometry.

\eject

\small
\vspace{1cm}
\tableofcontents{}
\vspace{1cm}
\bigskip\hrule
\normalsize

\section{Introduction}
\setcounter{equation}{0}
The moduli spaces of K\"ahler and complex structure deformations of
Calabi-Yau manifolds are naturally related to Special
K\"ahler (SK) and quaternion-K\"ahler (QK) geometry. Consequently,
these types of manifolds arise in the low energy effective action
of string theory compactifications on Calabi-Yau three-folds.

SK manifolds were discovered in the context of $N=2$ 
supergravity\footnote{The supersymmetry transformations have eight real components.}
theories coupled to vector multiplets \cite{deWit:1984pk}.
They are described by a holomorphic function
$F(X)$ that is homogeneous of degree two in complex coordinates
$X^I$. A more mathematically precise and intrinsic
formulation of this special geometry was given in
\cite{Freed:1997dp,Cortes:2001kx}. 

QK manifolds arose in the context of $N=2$
supergravity coupled to hypermultiplets \cite{Bagger:1983tt}.
QK manifolds are also described by a single function. This follows
from the construction of the Swann bundle \cite{Swann} over the QK
space. This bundle is hyperk\"ahler; locally, its metric is determined by a hyperk\"ahler
potential $\chi(\phi)$, where $\phi$ are local coordinates on the space. 
In \cite{DWRV}, we called such spaces hyperk\"ahler
cones (HKC's) because they have a homothety arising from
the underlying conformal symmetry. One therefore also
uses the terminology conformal hyperk\"ahler manifolds, as in
\cite{Bergshoeff:2004nf}.

In this note, we review the construction of Quaternion-K\"ahler
(QK) manifolds from special K\"ahler (SK) geometry, along the
lines of our recent work \cite{Rocek:2005ij}, but with more
emphasis on the mathematical structure. One constructs QK
spaces from SK manifolds by using the c-map \cite{CFG}. This
maps extends Calabi's construction of hyperk\"ahler metrics on
cotangent bundles of K\"ahler manifolds \cite{Calabi,CMMS}--a
construction well known in the mathematics community--to
quaternionic geometry.

As we shall see, the hyperk\"ahler cones arising from the c-map
have additional symmetries: they have an equal number of
commuting triholomorphic isometries as their quaternionic
dimension. Hyperk\"ahler manifolds with such isometries were
classified in \cite{Hitchin:1986ea} by performing a Legendre
transform on the (hyper)k\"ahler potential 
\cite{Lindstrom:1983rt} and writing the result
in terms of a contour integral of a meromorphic\footnote{Actually, certain branch
cut singularities sometimes arise.} function $H$. In
our case, because of the conformal symmetry of the HKC, 
this function $H$ is homogeneous of degree one. As a
result, the c-map induces a map from the holomorphic function $F$,
which characterizes the SK geometry and is
homogeneous of degree two, to a function $H$, which characterizes the QK geometry
and is homogeneous of degree one. Following
\cite{Rocek:2005ij}, we now describe this construction.

\section{The c-map}\label{cmap}
\setcounter{equation}{0}

In this section, we introduce our notation and review the c-map \cite{CFG,FS}.

Consider an affine (or rigid) special
K\"ahler manifold\footnote{We use the language of local
coordinates, but a coordinate free description can be found 
in \cite{Freed:1997dp,Cortes:2001kx}.} of dimension $2(n+1)$.  It is characterized
by a holomorphic prepotential $F(X^I)$, which is homogeneous of
degree two ($I=1,\cdots,n+1$). The K\"ahler potential and metric
of the rigid special geometry are given by
\begin{equation}\label{K-pot}
K(X,\bar X)=i({\bar X}^IF_I-X^I{\bar F}_I)\ ,~~
{\rm d}s^2=N_{IJ}\,{\rm d}X^I{\rm d}\bar X^J\ ,~~
N_{IJ}=i(F_{IJ}-{\bar F}_{IJ})\ ,
\end{equation}
where $F_I$ is the first derivative of $F$, {\it etc}.

The projective (or local) special K\"ahler geometry is then of real
dimension $2n$, with complex inhomogeneous coordinates
\begin{equation}\label{proj-coord}
Z^I=\frac{X^I}{X^1} = \{1,Z^A\}\ ,
\end{equation}
where $A$ runs over $n$ values. Its K\"ahler potential is given
by
\begin{equation}\label{local-K}
{\cal K}(Z,\bar Z)= {\rm ln} (Z^IN_{IJ}{\bar Z}^I)\ .
\end{equation}
We further introduce the matrices \cite{deWit:1984pk}
\begin{equation}\label{curlyN}
{\cal N}_{IJ}=-i{\bar F}_{IJ}- \frac{(NX)_I(NX)_J}{(XNX)}\ ,
\end{equation}
where $(NX)_I\equiv N_{IJ}X^J$, {\it etc}.

The c-map defines a $4(n+1)$-dimensional quaternion-K\"ahler
metric as follows: One builds a $G$-bundle, with $2(n+2)$ dimensional
fibers coordinatized by $\phi$, $\sigma$, $A^I$ and $B_I$, over the
projective special K\"ahler manifold; the real group $G$ is a semidirect product
of a Heisenberg group with $\mathbb{R}$, and acts on the coordinates by:
\begin{eqnarray}\label{G-action}
A^I\to {\rm e}^\beta(A^I+\epsilon^I)&,~~&
B_I\to {\rm e}^\beta(B_I +\epsilon_I)~,\nonumber\\[2mm] \phi\to\phi+\beta&,~&
\sigma\to {\rm e}^{2\beta}(\sigma+\alpha-\ft12 \epsilon_IA^I+\ft12\epsilon^IB_I)\ ,
\end{eqnarray}
Then the explicit $G$-invariant QK metric is \cite{FS}
\begin{eqnarray}\label{QK-metric}
{\rm d}s^2&=&{\rm d}\phi^2 - {\rm e}^{-\phi}({\cal N}+{\bar {\cal N}})_{IJ}
W^I{\bar W}^J+{\rm e}^{-2\phi}\Big( {\rm d}\sigma
-\frac{1}{2}(A^I{\rm d}B_I-B_I{\rm d}A^I)\Big)^2
\nonumber\\
&&-4 {\cal K}_{A\bar B}\,{\rm d}Z^A {\rm d}{\bar Z}^{\bar B}\ .
\end{eqnarray}
The metric is only positive definite in the domain where
$(ZN{\bar Z})$ is positive and hence ${\cal K}_{A{\bar B}}$ is negative
definite. One can then show that
${\cal N}+{\bar {\cal N}}$ is negative definite \cite{CKVPDFDWG}.  The
one-forms $W^I$ are defined by
\begin{equation}
W^I=({\cal N}+{\bar {\cal N}})^{-1\,IJ}\Big(2{\bar {\cal N}}_{JK}
{\rm d}A^K-i{\rm d}B_J\Big)\ .
\end{equation}
As shown in \cite{FS}, such metrics are indeed quaternion-K\"ahler;
they were further studied in \cite{DWVVP}, including an analysis of their isometries.
There are always the $2(n+2)$ manifest isometries (\ref{G-action}), 
of which $n+2$ are commuting, {\it e.g.,}
\begin{equation}\label{QK-isometries}
B_I\rightarrow B_I +\epsilon_I \ ,\qquad \sigma\rightarrow \sigma
+\alpha-\ft12 \epsilon_IA^I\ .
\end{equation}
This is one isometry {\em more} than the
quaternionic dimension of the QK manifold.

\section{Hyperk\"ahler cones and the Legendre transform}
\label{Superspace}
\setcounter{equation}{0}

The Swann bundle over a QK geometry, {\it i.e.,} the hyperk\"ahler cone
(HKC), is a hyperk\"ahler manifold with one extra quaternionic dimension. 
As for special K\"ahler manifolds, the geometry of the HKC is again affine. 
In physics terminology, this arises because
one adds a compensating hypermultiplet. Adding the compensator to
the original hypermultiplets that parametrize the
$4(n+1)$-dimensional QK space, one obtains a cone with real
dimension $4(n+2)$. This space is hyperk\"ahler and admits a
homothety as well as an $SU(2)$ isometry group that rotates the three
complex structures.

The metric on the HKC can be constructed from a hyperk\"ahler potential
\cite{Swann}, which is a K\"ahler potential with respect to {\em any} of
the complex structures. In real local coordinates $\phi^A$, the
metric and the hyperk\"ahler potential $\chi(\phi)$ are related by
\begin{equation}\label{chi}
g_{AB}=D_A \partial_B \,\chi(\phi)\ ,
\end{equation}
where $D_A$ is the Levi-Civita connection. As for all K\"ahler
manifolds, in complex coordinates, the hermitian part of (\ref{chi})
defines the metric in terms of the complex hessian of the potential; however,
in this case, the vanishing of the holomorphic parts of the metric
is an additional constraint on the geometry.

Any QK isometry can be lifted to a triholomorphic
isometry on the HKC. In
the physics literature, this was shown in
\cite{DWRV,deWit:2001bk}. Using the notation of the previous
section, we thus have an HKC of real dimension $4(n+2)$ together
with $n+2$ commuting triholomorphic isometries determined by
(\ref{QK-isometries}). As mentioned before, hyperk\"ahler
manifolds of this type were classified in \cite{Hitchin:1986ea}.
It is convenient to introduce complex coordinates $v^{\hat I}$ and
$w_{\hat I}$, in such a way that the isometries act as imaginary
shifts in $w_{\hat I}$. Notice that $\hat I=0,1,...,n+1$. The
hyperk\"ahler potential is then a function $\chi(v,\bar v,w+\bar
w)$, and can be written as a Legendre transform of a function
${\cal L}(v,\bar v,G)$ of $3(n+2)$ variables, The Legendre
transform with respect to $G^{\hat I}$ is defined by
\begin{equation}\label{legendre-tr}
\chi (v,\bar v,w,\bar w) \equiv {\cal L} (v,\bar v,G)  - (w+\bar
w)_{\hat I}\,G^{\hat I}\ , \qquad w_{\hat I}+\bar w_{\hat I} = 
\frac{\partial {\cal L}}{\partial G^{\hat I}}\ .
\end{equation}
The constraints from hyperk\"ahler geometry can be solved by
writing ${\cal L}$ in terms of a contour integral
\cite{PSS1,PSS2,Hitchin:1986ea}
\begin{equation}\label{c-int}
{\cal L}(v,\bar v,G)\equiv{\rm Im} \oint_{\cal C} \frac{{\rm d}
\zeta} {2 \pi i \zeta}\; H(\eta, \zeta)\ ,
\end{equation}
with
\begin{equation}\label{eta}
\eta^{\hat I} \equiv \frac{v^{\hat I}}{\zeta} + G^{\hat I} - {\bar v}^{\hat I}
\zeta\ .
\end{equation}
These objects have an interpretation in twistor space as sections
of an ${\cal O}(2)$ bundle. In physics terminology,
these are $N=2$ tensor multiplets.
Furthermore, the conditions for a homothetic Killing vector and
$SU(2)$ isometries imply that $H$ is a function homogeneous of
first degree\footnote{Actually, quasihomogeneity up to terms of
the form $\eta\ln(\eta)$ is sufficient \cite{DWRV}, but such terms
do not seem to arise in the c-map.}  (in $\eta$) and without
explicit $\zeta$ dependence \cite{DWRV}.

Since $H$ is homogeneous of first degree in $\eta$, it follows that the
hyperk\"ahler potential is also homogeneous of first degree in $v$ and
$\bar v$:
\begin{equation}\label{chi-homog}
\chi(\lambda v,\lambda {\bar v},w,\bar w)=\lambda \chi(v,{\bar v},w,\bar w)\ .
\end{equation}
In addition to a homothety, hyperk\"ahler cones also have an $SU(2)$ 
isometry group that rotates the sphere of complex structures. 
Under infinitesimal variations with respect to an element of the Lie algebra
$\varepsilon^+T_++\varepsilon^-T_-+\varepsilon^3T_3$, with 
$\varepsilon^-=(\varepsilon^+)^*$, 
these act as \cite{DWRV}
\begin{equation}\label{su2trans}
\delta_{\vec\varepsilon\,\,} v^{\hat I} =- i\varepsilon^3 v^{\hat I} +\varepsilon^- G^{\hat I}
(v,{\bar v},w,{\bar w})\ , \qquad \delta_{\vec\varepsilon\,\,} w_{\hat I}=\varepsilon^+ 
\frac{\partial \cal L}{\partial {\bar v}^{\hat I}}\ ,
\end{equation}
where $G^{\hat I}$ has to be understood as the function of the coordinates
$v,{\bar v},w,{\bar w}$ obtained by the Legendre transform defined in
(\ref{legendre-tr}). The coordinates $w_I$ do not transform under
$\varepsilon^3$. One can now explicitly check that the hyperk\"ahler
potential is $SU(2)_R$ invariant,
\begin{equation}
\delta_{\vec\varepsilon\,\,} \chi = {\cal L}_{v^{\hat I}}\,\delta_{\vec\varepsilon\,\,} v^{\hat I}
+{\cal L}_{{\bar v}^{\hat I}}\,\delta_{\vec\varepsilon\,\,} {\bar v}^{\hat I}
-\delta_{\vec\varepsilon\,\,} (w_{\hat I}+{\bar w}_{\hat I})\, G^{\hat I}=0\ .
\end{equation}
(The $\delta G$ terms cancel identically because $\chi$ is a Legendre 
transform). For the generators $\varepsilon^{\pm}$ this is immediately 
obvious; for variations proportional to $\varepsilon^{3}$ one needs to use the
invariance of ${\cal L}$, \ie $v^{\hat I}{\cal L}_{v^{\hat I}} = 
{\bar v}^{\hat I}{\cal L}_{{\bar v}^{\hat I}}$.

\section{Hyperk\"ahler cones from the c-map}

The Quaternion-K\"ahler space in the image of the c-map has
dimension $4(n+1)$. The hyperk\"ahler cone above it has dimension $4(n+2)$.
It therefore needs to be described by $n+2$ twistor variables, say 
$\eta^I$ and $\eta^0$, where $I=1,\cdots , n+1$.
As we shall show, the result for the tree level c-map is
given by
\begin{equation}\label{G-cmap}
H(\eta^I,\eta^0) = \frac{F (\eta^I)}{\eta^0}\ ,
\end{equation}
where $F$ is the prepotential of the special K\"ahler geometry, now evaluated
on the twistor variables $\eta$. This is our main result.
Note that $H$ does not depend explicitly on $\zeta$ and, since $F$ is 
homogeneous of degree two, $H$ is homogeneous of degree one, as required by 
superconformal invariance.

We now give a detailed proof of (\ref{G-cmap}) by explicit
calculation \cite{Rocek:2005ij}. To be precise, we prove that
(\ref{G-cmap}) leads to (\ref{QK-metric}).

\subsection{Gauge fixing and the contour integral}

As explained in the previous section, any hyperk\"ahler cone has an
$SU(2)$ symmetry and a homothety. The generators of the homothety and
$U(1)\subset SU(2)$ give a natural complexified action on the HKC;  the
remaining two generators of the $SU(2)$ combine to give the roots $T_{\pm}$.

To evaluate the contour integral (\ref{c-int}), it is convenient to 
make use of the isometries. In physics terminology, one 
can impose gauge choices. Mathematically, the isometries fiber the total
space by the orbits, and a gauge choice is just a choice of section. 
For the  symmetries generated by $T_{\pm}$, whose action is 
given by (\ref{su2trans}), we choose 
\begin{equation}
v^0 = 0 \ .
\end{equation}
In this gauge, we have that $\eta^0=G^0$
and this simplifies the pole structure in the complex
$\zeta$-plane. Then, using as well the homogeneity properties of $F$,
the contour integral (\ref{c-int}) simplifies to 
\begin{equation}
{\cal L} (v,{\bar v},G)= \frac{1}{G^0}\, {\rm Im}\ \oint
\frac{{\rm d} \zeta}{2\pi i} \frac{F(\zeta \eta^I)}{\zeta^3}\ ,
\end{equation}
with
\begin{equation}
\zeta \eta^I={v}^I+\zeta G^I - \zeta^2 {\bar v}^I\ ,\qquad
I=1,\cdots , n+1\ ,
\end{equation}
which, for nonzero values of $v$, has no zeroes at $\zeta=0$.
Therefore, since $F(\zeta \eta)$ is homogeneous of positive degree, it
has no poles at $\zeta=0$.
It is now easy to evaluate the contour integral, because
the residue at $\zeta=0$ replaces all the $\zeta \eta^I$ by $v^I$.
The result is
\begin{equation}\label{cont-F}
{\cal L}(v,\bar v,G)=\frac{1}{4G^0}\Big( N_{IJ}G^I G^J - 2 K(v,\bar v) \Big)\ ,
\end{equation}
where $K(v,\bar v)$ is the K\"ahler potential of the rigid special
geometry given in (\ref{K-pot}), with $F_I(v)$ now the derivative
with respect to $v^I$, {\it etc}. Notice that the function ${\cal L}$ satisfies the 
Laplace-like equations \cite{PSS1,PSS2,Hitchin:1986ea}
\begin{equation}
{\cal L}_{G^IG^J}+{\cal L}_{v^I{\bar v}^J}=0\ .
\end{equation}

The equation is not satified for the components ${\cal L}_{G^0G^0}$
and ${\cal L}_{G^0G^I}$, because we have chosen the gauge
$v^0=0$. It would be interesting to compute ${\cal L}$ for arbitrary
values of $v^0$. For a special case, this was done in \cite{Anguelova:2004sj}.

\subsection{The hyperk\"ahler potential}

To compute the hyperk\"ahler potential $\chi$, we have to 
Legendre transform ${\cal L}$,
\begin{equation}\label{legendre-chi}
\chi (v,\bar v,w,\bar w) = {\cal L}(v,\bar v,G) +(w+{\bar w})_0\,G^0  -
(w+{\bar w})_I\,G^I \ ,
\end{equation}
The hyperk\"ahler potential $\chi$, computed by 
extremizing\footnote{The relative minus signs between the
last two terms in (\ref{legendre-chi}) is purely a matter of
convention.} (\ref{legendre-chi}) with respect  $G^0,G^I$
completely determines the associated hyperk\"ahler geometry. In general, it 
is a function of the $2(n+2)$
complex coordinates $v^0,v^I$ and $w_0,w_I$, but we 
work only on the (K\"ahler but not hyperk\"ahler) submanifold $v^0=0$.
This is sufficient for calculating the metric on the underlying QK manifold.
The geometry of the HKC only depends on $w$ through the
combination $w+\bar w$ which makes manifest
the $n+2$ commuting isometries. The Legendre transform of (\ref{cont-F})
gives:
\begin{equation}\label{LT-phi}
\frac{G^I}{G^0}=2N^{IJ}(w+{\bar w})_J\ ,\qquad (G^0)^2=
\frac{K}{2\Bigl((w+{\bar w})_IN^{IJ}(w+\bar w)_J-(w+{\bar w})_0\Bigr)}\ .
\end{equation}
Up to an irrelevant overall sign, we find, using (\ref{cont-F})
\begin{equation}\label{chi-phi}
\chi \Bigl(v,\bar v, G(v,\bar v,w,\bar w)\Bigr)=\frac{K(v,\bar
v)}{G^0}\ ,
\end{equation}
where $G^0$ is determined by (\ref{LT-phi}). More explicitly,
in terms of the HKC coordinates,
\begin{equation}\label{cmap-HKC}
\chi (v,\bar v, w,\bar w)= {\sqrt 2}\,\,{\sqrt {K(v,\bar v)}
}\,{\sqrt {(w+{\bar w})_IN^{IJ}(w+\bar w)_J-(w+{\bar w})_0}}\ .
\end{equation}

\subsection{Twistor space}
The twistor space above a $4(n+1)$ dimensional QK manifold has dimension
two higher, and is K\"ahler. It can be seen as a CP$^1$ bundle
over the QK. It can also be obtained from the HKC by a K\"ahler quotient
with respect to $U(1)\subset SU(2)$.
Equivalently, we define inhomogeneous coordinates, {\it e.g.},
\begin{equation}
Z^I=\frac{v^I}{v^1}=\{1,Z^A\}\ ,
\end{equation}
where $A$ runs over $n$ values. As we show below, these
inhomogeneous coordinates will be identified with
(\ref{proj-coord}).

The K\"ahler potential on the twistor space, denoted by $K_T$, is
given by the logarithm of the hyperk\"ahler potential restricted to the submanifold
given by $v^1=1$ \cite{DWRV}:
\begin{equation}\label{ktzw}
K_T(Z,\bar Z,w,\bar w)=\frac{1}{2}\Big[ {\cal K}(Z,\bar Z)
+ {\rm ln} \Big((w+{\bar
w})_IN^{IJ}(w+\bar w)_J-(w+{\bar w})_0\Big) \Big] +{\rm ln}({\sqrt
2})\ ,
\end{equation}
where ${\cal K}(Z,\bar Z)$ is the {\em same} as the special K\"ahler potential
(\ref{local-K}).

On the twistor space, there always exists a holomorphic one-form
${\cal X}$ which can be constructed from the holomorphic two-form
that any hyperk\"ahler manifold admits. In our case this one-form
is obtained from the holomorphic HKC two-form $\Omega = {\rm d}w_I
\wedge {\rm d}v^I$. Without going into details, it is given by
\cite{DWRV}
\begin{equation}
{\cal X}=2Z^I{\rm d}w_I\equiv {\cal X}_\alpha {\rm d}z^\alpha\ ,
\end{equation}
where the index $\alpha=1,\cdots ,2(n+1)$ runs over the complete
set of holomorphic coordinates $w_I,w_0,Z^A$ on the submanifold\footnote{This
submanifold can be thought of invariantly as a quotient of the original HKC.}
of the twistor space given by $v^0=0$. In total this gives 
$2(n+1) + 2 + 2n = 4(n+1)$--the (real) dimension of the QK.
The metric on the QK manifold can then be 
computed\footnote{Note that the constant
term in $K_T$ (\ref{ktzw}) enters in (\ref{G-metric}).}:
\begin{equation}\label{G-metric}
G_{\alpha{\bar \beta}}=K_{T,\,\alpha {\bar \beta}}-{\rm
e}^{-2K_T}{\cal X}_\alpha {\bar {\cal X}}_{{\bar \beta}}\ .
\end{equation}

\subsection{The quaternionic metric}

We now compute the QK metric that follows from the
c-map using (\ref{G-metric}). To compare with (\ref{QK-metric})
we only need to identify the coordinates $w_I,w_0$ with those of
(\ref{QK-metric}), since the $Z^A$ coordinates of the special
K\"ahler manifold can be identified with the ones above. We define
\begin{eqnarray}\label{QK-coord1}
w_0&=&iA^IA^JF_{IJ}-i(\sigma+\frac12A^IB_I)-{\rm e}^\phi\ ,\nonumber\\
w_I&=&iF_{IJ}A^J-\frac{i}{2}B_I\ .
\end{eqnarray}
The metric can be written in these coordinates; after considerable calculation
\cite{Rocek:2005ij}, up to an irrelevant overall all factor of $-1/8$, 
we obtain precisely the result (\ref{QK-metric})!
From (\ref{LT-phi}), we find the following
relations between the QK coordinates and the twistor variables $G^{\hat I}$ (see \ref{eta}):
\begin{equation}\label{expr-dilaton}
2A^I =\frac{G^I}{G^0}\ ,\qquad 4{\rm e}^\phi = \frac{K(v,\bar v)}{(G^0)^2}\ .
\end{equation}
This concludes the proof of (\ref{G-cmap}).

\section{Summary and conclusion}
\setcounter{equation}{0}

We have constructed the Swann bundle over the Quaternion-K\"ahler manifolds
that arise in the c-map. The corresponding hyperk\"ahler 
potential was given in (\ref{cmap-HKC}), and was first derived in 
\cite{Rocek:2005ij}. 
Introducing coordinates
\begin{equation}\label{new-X}
X^I(v,\bar v,w,\bar w)\equiv \frac{v^I}{{\sqrt {G^0(v,\bar v,w,\bar w)}}}\ ,
\end{equation}
we can conveniently rewrite the hyperk\"ahler potential as
\begin{equation}\label{chi-is-K}
\framebox{$\chi(v,\bar v,w,\bar w)=K\Big(X^I(v,\bar v,w,\bar w),{\bar X}^I(v,\bar v,w,\bar w)\Big)$}\ .
\end{equation}
Here $K$ is the K\"ahler potential $K(X,\bar X)=i({\bar X}^IF_I-X^I{\bar F}_I)$ of the affine special geometry.

The special hyperk\"ahler cones given by the c-map have as 
many ($n+2$) commuting triholomorphic
isometries as their quaternionic dimension. As explained before, this implies
the hyperk\"ahler potential can be Legendre transformed to a 
function ${\cal L}$ that can be written in terms of a contour integral over 
a function $H(\eta)$; equivalently, the twistor space of the HKC can be described
in terms of sections of $(n+2)$ ${\cal O}(2)$-bundles. These twistor variables 
$\eta$ were defined in (\ref{eta}) and the function 
$H$ was determined in (\ref{G-cmap}). Defining
\begin{equation}
X^I(\eta) \equiv \frac{\eta^I}{{\sqrt {\eta^0}}}\ ,
\end{equation}
we can write $H$ as
\begin{equation}\label{H-is-F}
\framebox{$H(\eta^I,\eta^0)=F\big[X^I(\eta^I,\eta^0)\big]$}\ .
\end{equation}

The function $F$ is well known to be related to the topological 
string amplitude \cite{Antoniadis:1993ze,Bershadsky:1993cx}. Typical
examples that appear in the context of Calabi-Yau compactifications 
are of the form
\begin{equation}
F(X^I)=d_{ABC}\frac{X^AX^BX^C}{X^1}\ ,
\end{equation}
where $X^I=\{X^1,X^A\}$ and the constants $d_{ABC}$ are related to the
triple intersection numbers of the Calabi-Yau. To give an 
explicit example, one can choose specific
values for these coefficients such that the local (projective) special 
K\"ahler geometry is  the symmetric space
\begin{equation}
\frac{SU(1,1)}{U(1)}\times \frac{SO(n-1,2)}{SO(n-1)\times SO(2)}\ .
\end{equation}
After the c-map, the hyperk\"ahler cone is based on
the function
\begin{equation}
H(\eta^{\hat I})=d_{ABC}\frac{\eta^A\eta^B\eta^C}{\eta^0\eta^1}\ .
\end{equation}
This corresponds to a homogeneous quaternion-K\"ahler manifold of 
the form (see e.g. appendix C in \cite{CFG}, and references therein)
\begin{equation}
\frac{SO(n+1,4)}{SO(n+1)\times SO(4)}\ .
\end{equation}
Other examples were recently given in \cite{RSV}, were quantum effects
were taken into account.

The connection of these geometries with topological strings is very profound,
and has important physical implications. For instance, it was 
recently shown that the topological string amplitude $F$ appears
in the study of supersymmetric black holes in string 
theory \cite{Ooguri:2004zv,OVV}. More precisely, the Legendre transform
of $F$ is related to the entropy of the black hole. It would be 
interesting to see if this Legendre transform is related to the one described
here; speculations along these lines can be found in \cite{Rocek:2005ij}. 
To make progress on this issue, one needs to evaluate the contour integral (\ref{c-int}) 
without making use of the special coordinate system in which we can 
set $v^0=0$. We leave this for future research.

\vskip .5cm
\noindent{\bf Acknowledgements}

\noindent 

This work was presented at the 77th ``Rencontre entre Physiciens
Th\'eoriciens et Math\'ematiciens" on ``Pseudo-Riemannian Geometry
and   Supersymmetry", Strasbourg. SV thanks V.~Cort\'es for the 
invitation and kind hospitality. 

Most of this work has been initiated and completed during the
2004 and 2005 Simons Workshops in Physics and Mathematics.
SV and CV thank the C.N. Yang Institute for Theoretical Physics
and the Department of Mathematics at Stony Brook University for hosting the
workshops and for partial support. MR thanks the Institute for Theoretical Physics  
at the University of Amsterdam for hospitality. MR is supported in part by
NSF grant no.~PHY-0354776, by the University of Amsterdam, and by Stichting FOM. 
CV is supported in part by NSF grants PHY-0244821 and DMS-0244464.

\raggedright

\end{document}